\documentclass[a4paper,11.5pt,reqno]{amsart}

\usepackage{graphicx}
\usepackage[utf8]{inputenc}
\usepackage{amsmath,amssymb,amsthm}
\usepackage{enumitem}
\usepackage{xcolor}
\usepackage{hyperref}

\usepackage{lineno}
\modulolinenumbers[3]

\theoremstyle{plain}
\newtheorem{theorem}{Theorem}[section]

\theoremstyle{definition}
\newtheorem{definition}{Definition}

\theoremstyle{remark}
\newtheorem{remark}{Remark}

\newcommand{\pMIo}{\mathbb{M}}

\numberwithin{equation}{section}

\hypersetup{colorlinks = true, linkcolor = blue, urlcolor = blue, citecolor = blue}

\begin{document}

\title[On two Abelian Groups Related to the Galois Top] 
      {On two Abelian Groups Related to the \\ Galois Top}

\author[H. Ruhland]{Helmut Ruhland}
\address{Santa F\'{e}, La Habana, Cuba}
\email{helmut.ruhland50@web.de}



\begin{abstract}
In mathematical physics the Galois top, introduced by S. Adlaj, possesses a fixed point on one of two Galois axes through its center of mass. This heavy top has two algebraic motion invariants and an additional transcendental motion invariant. This third invariant depends on an antiderivative of a variable in the canonical phase space.  \\
In this article an abelian semigroup and an abelian group are defined that are related to the application of the Huygens-Steiner theorem to points on the Galois axis of a rigid body. This yields nonlinear representations of the one-dimensional, affine, linear (semi)group.
\end{abstract}

\date{\today}

\maketitle

\section{Introduction}

The author of \cite{Adlaj} introduced the Galois top. The fixed point $O$ of this top is located on one of the two so-called Galois axes through the center of mass $G$ of the rigid body. A Galois axis can be interpreted geometrically in this manner: it passes through $G$ and is orthogonal to a plane that intersects the MacCullagh ellipsiod centered at $G$ in a circle. Besides the two motion invariants - energy $K$ and angular momentum projected onto the vector of gravity $L_g$ (every heavy top has these two invariants) - the Galois top has a third \textit{transcendental} constant of motion. It is widely \textit{believed} that such a transcendental invariant (depending on an antiderivative of a variable in the canonical phase space) \textit{would not exist}; for the invariant see formula (3) in \cite{Adlaj}. \\ 

The application of the Huygens-Steiner theorem to points $O$ on the Galois axis defines maps from the principal moments of inertia in the center of mass $G$ to the principal moments of inertia in the
point $O$. \\
In the next section, it is shown that these maps generate an abelian semigroup. In section \ref{section_group} an abelian group is described, also generated by these maps. Finally in section \ref{section_group_ext} nonabelian extensions of this semigroup and group are given. \\

To give a short summary of the findings in this article in one sentence: \\
Applying the Huygens-Steiner theorem in a physical setting, we construct a \textit{nonlinear} representation of the one-dimensional, affine, linear (semi)group, the so-called $a x + b$ group.
 
\section{An abelian semigroup of inertia maps \label{section_semigroup}}

Let $\pMIo \subset \mathbb{R}^3$ represent the three (ordered with respect to $<$) principal moments of inertia of a rigid body, $\pMIo = \{ ( A, B, C ) \in \mathbb{R}^3 \, \vert \, 0 < A < B < C \}$. The condition $0 < A$ in the definition ensures that the 3 principal moments are physically meaningful.  \\

Then the rwo-valued (2 Galois axes exist) tensor of inertia $J(d)$, defined in the line below formula (2) in \cite{Adlaj}, defines maps acting on $\pMIo$:

\begin{definition}
For $x \in \mathbb{R}_{\geq 0}$ this is a one-parameter family of maps:
\begin{align}
 j(x) : \enspace  & \hspace{0.6cm} \pMIo & \rightarrow & \hspace{0.8cm} \pMIo \nonumber \\
 & ( A, B, C ) & \mapsto & \enspace ( \lambda_1, \lambda_2, \lambda_3 )  \label{jmap_lambda}
\end{align}
 $0 < \lambda_1 < \lambda_2 = B + x < \lambda_3$ are the three ordered eigenvalues of $J(d)$, $d^2$
 replaced by $x$. These maps depend only on $x$, not on which of the two Galois axes is chosen. \\
 In Appendix \ref{proof_codomain}, it is shown that the condition $x \ge 0$ is necessary and sufficient.
\end{definition}

\begin{theorem} \label{theorem_j}
The maps $j(d^2)$, where $d$ is the distance of a point $O$ on the Galois axis to the center of mass $G$,
define an abelian semigroup $G_s = \{ j(d^2) \, \mid \, d \in  \mathbb{R} \} =
 \{ j(x) \, \vert \, x \in \mathbb{R}_{\geq 0} \}$. $j(0)$ is the neutral element. The semigroup law is
$j(x) \circ j(y) = j(x + y)$. For a proof, see Appendix \ref{proof_semigroup}.
\end{theorem}

\section{An abelian group of maps \label{section_group}}

Starting with the semigroup $G_s$ defined in Theorem \ref{theorem_j}, we could try to get an abelian group $G = \{ j(x) \, \vert \, x \in  \mathbb{R} \}$ by adding the inverses of all elements in $G_s$. The inverse of $j(x)$ would be $j(- x)$. But $j(x)$ is only defined for $x \ge 0$ and not for $x \in \mathbb{R}$, as  would be necessary to get a group $G$. \\

But we can define other maps with $\mathbb{C}^3$ as domain and codomain. Based on these maps, we obtain an abelian group. In contrast to the semigroup $G_s$ above, this group no longer has a reference to physical tops. \\

The tensor of inertia $J(d)$ now defines maps acting on $\mathbb{C}^3$:

\begin{definition}
For $x \in \mathbb{C}$ this is a one-parameter family of maps:
\begin{align}
 \overline{j}(x) : \enspace  & \hspace{0.6cm} \mathbb{C}^3 & \rightarrow & \hspace{1.0cm} \mathbb{C}^3 \nonumber \\
 & ( A, B, C ) & \mapsto & \enspace ( \lambda_{1/3}, \lambda_2, \lambda_{3/1} )  \label{jbar_map_lambda}
\end{align}
$\lambda_1, \lambda_2 = B + x, \lambda_3$ are the three eigenvalues of $J(d)$, $d^2$ replaced by $x$. \\
These maps are two-valued/sheeted after analytic continuation. An involution $i$ leaves such a two-valued map invariant; it interchanges the two sheets or interchanges the components $A, C$ of the domain:
\begin{equation}
i : \mathbb{C}^3 \rightarrow \mathbb{C}^3, \; ( A, B, C ) \mapsto ( C, B, A )
\end{equation} \\
$i \circ \overline{j}(x) = \overline{j}(x)$ the 2 sheets are interchanged, i.e., $( \lambda_3, \lambda_2, \lambda_1 ) = i \, ( \lambda_1, \lambda_2, \lambda_3 )$.
$\overline{j}(x) \circ i = \overline{j}(x)$: $A, B$ are interchanged, i.e., $( C, B, A ) = i \, ( A, B, C )$.
\end{definition}

\begin{theorem}
The maps $\overline{j}(x)$ define an abelian group $G = \{ \overline{j}(x) \, \mid \, x \in  \mathbb{C} \}$.
$\overline{j}(0)$ is the neutral element. The inverse of $\overline{j}(x)$ is $\overline{j}(-x)$. The group law is $\overline{j}(x) \circ \overline{j}(y) = \overline{j}(x + y)$.
\end{theorem}

\section{A (semi)group extension by the scaling group \label{section_group_ext}}

\begin{definition}
For $x \in \mathbb{R_+}/\mathbb{C}^\times$, with $\mathbb{C}^\times = \mathbb{C} \setminus \{ 0 \}$,  the following maps build a one-parameter family, denoted as scaling maps. The alternative in the domain/codomain depends on the context:
\begin{align}
 s(x)/\overline{s}(x) : \enspace  & \hspace{0.3cm} \mathbb{M}/\mathbb{C}^3 & \rightarrow & \hspace{0.6cm} \mathbb{M}/\mathbb{C}^3 \nonumber \\
 & ( A, B, C ) & \mapsto & \enspace ( x \, A, x \, B, x \, C )  \label{s_map}
\end{align}
\end{definition}
The abelian group $S = \{ s(x)/\overline{s}(x) \, \vert \, x \in  \mathbb{R_+}/\mathbb{C}^\times \}$ is denoted as scaling group. \\

The inverse of $s(x)$ is $s(1/x)$. We have the following identity:
\begin{equation}
s(x_1) \circ j(x_2) \circ s(x_1)^{-1} = s(x_1) \circ j(x_2) \circ s(1 / x_1) = j(x_2 / x_1)
\end{equation}
Therefore we can extend the semigroup $G_s$ and the group $G$ to the two-parameter nonabelian semigroup $G_s.S$ and nonabelian group $G.S$. Using left cosets: $G_s.S = \{ s(x_1) \circ j(x_2) \, \vert \, x_1 \in  \mathbb{R_+} \, \wedge \, x_2 \in \mathbb{R}_{\geq 0} \}$ and  $G.S = \{ s(x_1) \circ \overline{j}(x_2) \, \vert \, x_1 \in  \mathbb{C}^\times \, \wedge \, x_2 \in \mathbb{C}  \}$. \\

The group multiplication is now:
\begin{equation}
( s(x_1) \circ j(y_1) ) \quad \circ \quad ( s(x_2) \circ j(y_2) ) = ( s(x_1 x_2) \circ j(y_1 + y_2 \, x_1) )
\label{group_mult}
\end{equation} \\

\begin{remark}
These (semi)groups $G_s.S$ and $G.S$ are isomorphic to the two-parameter affine, linear (semi)group defined by the maps $al(a, b) : x \mapsto a \, x + b$. It can be verified easily that the latter group has the same (semi)group multiplication as (\ref{group_mult}).
\end{remark}

\section{An open question}

Perhaps it can be shown that arbitrary axes through the center of mass $G$, except the two Galois axes,
do not allow the definition of such abelian (semi)groups. \\

Besides the characterization of the two Galois axes using sections of the MacCullagh ellipsiod, this would allow us to characterize the Galois axes as the only axes with assigned (semi)groups.

\newpage
\noindent \textbf{\large Appendices}

\appendix

\section{Proof that the image of $j(x)$ is a subset of the codomain \label{proof_codomain}}

It must be proven: $\{ j(x) \, m \, \vert \, m \in \pMIo \} \subseteq \pMIo$ for $x \ge 0$. \\

Let $K(x)$ be the symmetric $2 \times 2$ submatrix of $J(d)$ with rows/columns 1 and 3, where $d^2$ is replaced by $x$:
\begin{equation}
   K (x) = 
   \begin{pmatrix}
     A + \frac{A \, (B - C) \, x}{B \, (A - C)}  &
     \frac{\sqrt{C A \, (A - B) \, (B - C)} \, x}{B \, (A - C)} \\
     \frac{\sqrt{C A \, (A - B) \, (B - C)} \, x}{B \, (A - C)} &
     C + \frac{C \, (A - B) \, x}{B \, (A - C)}
   \end{pmatrix}
\end{equation}
Then:
\begin{equation}
 Tr_x = \mathrm{trace}(K(x)) = A + C + x,  \quad Det_x = \det(K(x)) = A \, C \, (1 + x / B)
\end{equation}
The characteristic equation of $K(x)$ is then: $\lambda ^2 - Tr_x \, \lambda + Det_x = 0$. \\

The two eigenvalues of $K(x)$, and thus the lowest and highest eigenvalues of $J(d)$, are:
\begin{equation}
 \lambda_{1/3} = ( Tr_x   \mp   \sqrt{\Delta_x} ) / 2 \qquad \Delta_x = Tr_x ^ 2 - 4 \, Det_x 
\end{equation}
For the square root $\sqrt{\Delta_x}$, the principal value has to be taken. \\

Map \ref{jmap_lambda} now has this form:
\begin{align}
j(x) : \enspace  & \hspace{0.6cm} \pMIo & \rightarrow & \hspace{4.0cm} \pMIo \nonumber \\
& ( A, B, C ) & \mapsto & \enspace ( A_x = (Tr_x - \sqrt{\Delta_x}) / 2,     B_x = B + x,    C_x = (Tr_x + \sqrt{\Delta_x}) / 2 ) \label{jmap_detailed}
\end{align}
Because $A_x$ and $C_x$ are eigenvalues of $K(x)$, see the characteristic equation above, we have: $A_x + C_x = Tr_x, \; A_x \, C_x = Det_x$. The determinant of the Jacobian $\det\mathrm{jac}(j(x))$ is
$(C - A) \, (1 + x / B) / \sqrt{\Delta_x}$. \\

We prove the following inequalities for $x \ge 0$:

\begin{enumerate}

\item The square root $\sqrt{\Delta_x}$ defined above is real, i.e., $\Delta_x$ is positive: \\
$\Delta_x = x ^ 2 + ...$ is a quadratic polynomial in $x$. Solving $d\Delta_x / dx = 0$ gives $x_{min}$;
the minimum is  $4 \, A \, C \, (B - A) \, (C - B) / B ^ 2$ and is positive\footnote{a "physical proof" that $\Delta_x$ is positive for $x \ge 0$: because the principal moments of inertia $A, B, C > 0$ in the center of mass are physical, the principal moments in all points are physical, i.e., real. This implies $\Delta_x$ is positive.} here for $x \in \mathbb{R}$.

\item $A_x < B_x$: \\
$(Tr_x - \sqrt{\Delta_x}) / 2 - B_x < 0$. With $D = A + B + x - 2 \, (B + x)$ we get $D < + \sqrt{\Delta_x}$.
$D ^ 2  - \Delta_x = - 4 \, (B + x) \, (B - A) \, (C - B) / B < 0$, hence $\vert D \vert < \sqrt{\Delta_x}$.
Therefore, for each of the two possible signs of $D$, $D < + \sqrt{\Delta_x}$, which implies $A_x < B_x$.

\item $B_x < C_x$: \\
$(Tr_x + \sqrt{\Delta_x}) / 2 - B_x > 0$. With $D$ as in (2), we get $D > - \sqrt{\Delta_x}$. \\
From (2), it was already shown that $\vert D \vert < \sqrt{\Delta_x}$,
so for each of the two possible signs of $D$, we have $D > - \sqrt{\Delta_x}$, which implies $B_x < C_x$.

\item $0 < A_x$: \\
In (3), it is already proven that $0 < B_x < C_x$. Moreover, $A_x \, C_x = Det_x = A \, C \, (1 + x / B) > 0$, which implies $0 < A_x$.
\end{enumerate}

Thus, $0 < A_x < B_x < C_x$, so the image lies in the codomain.  \qed

\section{Proof of the (semi)group law \label{proof_semigroup}}

It has to be proven: the composition of two maps $j(x), j(y)$ satisfies the group law $j(x) \circ j(y) = j(x + y)$ or $j(x) \, j(y) \, m = j(x + y) \, m$ for $m \in \pMIo$. \\

Using the formula (\ref{jmap_detailed}) for $j(x)$, the composition of two maps $j(x) \, j(y) \, (A, B, C)$ yields the vector:
\begin{equation} \begin{split}
( A_{xy}, B_{xy}, C_{xy}) \qquad & A_{xy} = ( A_y + B_y + x  - \sqrt{\Delta_{xy}} ) / 2, B_{xy} = B_y + x \\
                                 & C_{xy} = ( A_y + B_y + x  + \sqrt{\Delta_{xy}} ) / 2  \qquad \qquad A_y + C_y = Tr_y \\
                                 & A_{xy} = ( Tr_y + x  - \sqrt{\Delta_{xy}} ) / 2, B_{xy} = B + x + y  \\
                                 & C_{xy} = ( Tr_y + x + \sqrt{\Delta_{xy}} ) / 2  \qquad \qquad \qquad Tr_y + x = Tr_{x + y} \\
                                 & A_{xy} = ( Tr_{x + y} - \sqrt{\Delta_{xy}} ) / 2, B_{xy} = B + x + y \\
                                 & C_{xy} = ( Tr_{x + y} + \sqrt{\Delta_{xy}} ) / 2
\end{split} \end{equation}
 
The first indications of commutativity and additivity are already visible; it remains to calculate $\Delta_{xy}$:
\begin{equation} \begin{split}
 \Delta_{xy} & = (A_y + B_y + x) ^ 2   -  4 \, A_y \, C_y \, (1 + x / B_y)
																			 \qquad \qquad \qquad A_y \, C_y = Det_y  \\
             & = (Tr_y + x) ^ 2        -  4 \, A \, C \, (1 + y / B) \, (1 + x / (B + y)) \\
             & = Tr_{x + y} ^ 2        -  4 \, A \, C \, (1 + (x + y) / B) \\
             & = Tr_{x + y} ^ 2        -  4 \, Det_{x + y} = \Delta_{x + y}
\end{split} \end{equation}
 
Hence, $j(x) \, j(y) \, (A, B, C)  = j(x + y) \, (A, B, C)$.  \qed

\bibliographystyle{amsplain}

\end{document}